\documentclass{amsart}
\usepackage{amssymb,latexsym,enumerate}
\theoremstyle{plain}
\newtheorem{Thm}{Theorem}[section]

\newtheorem{Prop}[Thm]{Proposition}
\newtheorem{Lem}[Thm]{Lemma}

\newtheorem{Conj}[Thm]{Conjecture}
\numberwithin{equation}{section}
\DeclareMathOperator{\Mat}{M}

\DeclareMathOperator{\supp}{supp}

\begin{document}

\title[Von Neumann Algebras and Independence of Translates]
{Von Neumann Algebras and Linear Independence of Translates}

\author[P. A. Linnell]{Peter A. Linnell}

\address{Math \\ VPI \\ Blacksburg \\ VA 24061--0123
\\ USA}

\email{linnell@math.vt.edu}
\urladdr{http://www.math.vt.edu/people/linnell/}

\begin{abstract}
For $x,y \in \mathbb {R}$ and $f \in L^2(\mathbb {R})$, define
$(x,y) f(t) = e^{2\pi iyt} f(t+x)$ and
if $\Lambda \subseteq \mathbb
{R}^2$, define $S(f, \Lambda) = \{(x,y)f \mid (x,y) \in \Lambda \}$.
It has been conjectured that
if $f\ne 0$, then $S(f,\Lambda)$ is linearly
independent over $\mathbb {C}$; one motivation for this problem comes
from Gabor analysis.  We shall prove that $S(f, \Lambda)$ is linearly
independent if $f \ne 0$ and $\Lambda$ is contained in a discrete
subgroup of $\mathbb {R}^2$, and as a byproduct we shall obtain some
results on the group von Neumann algebra generated by the operators
$\{(x,y) \mid (x,y) \in \Lambda \}$.  Also we shall prove these
results for the obvious generalization to $\mathbb {R}^n$.
\end{abstract}

\keywords{group von Neumann algebra, Gabor analysis, Heisenberg
group}

\subjclass{Primary: 46L10; Secondary 42C99}

\maketitle

\section{Introduction} \label{Sintroduction}

Let $n$ be a positive integer, let
$\mathcal {G}_n$ be the abelian group
$\{(x,y) \mid x,y \in \mathbb {R}^n \}$ with the
operation addition (so $\mathcal{G}_n \cong \mathbb {R}^{2n}$),
and for $x,y \in \mathbb {R}^n$, let $x \cdot y$ denote the dot
product $x_1y_1 + \dots + x_n  y_n$.  Let
$\mathbb {C} *\mathcal {G}_n$
denote the twisted group ring (a twisted group ring is a particular
kind of crossed product) which has $\mathbb {C}$-basis $\{\bar g \mid
g\in \mathcal {G}_n\}$, and multiplication
satisfying $\overline{(a,b)}
\, \overline{(x,y)} = e^{2\pi ia \cdot y}\overline{(a+x, b+y)}$.
For $g\in \mathcal {G}_n$, we shall often write $g$ instead
of $\bar g$ if there is
no danger of confusion, and then $g^{-1}$ will mean $\bar g^{-1}$
rather than $\overline {g^{-1}}$.  Let $L^2(\mathbb
{R}^n)$ denote the Hilbert space of square integrable
functions $\{f \colon \mathbb {R}^n \to \mathbb {C} \mid
\int_{\mathbb {R}^n} |f(t)|^2 \,dt < \infty\}$ with two functions
$f_1,f_2 \in L^2(\mathbb {R}^n)$ being equal if and only if $f_1(t) =
f_2(t)$ almost everywhere, and let $\mathcal
B(L^2(\mathbb {R}^n))$ denote the set of bounded linear operators on
$L^2(\mathbb {R}^n)$.  Then $\mathbb {C}*\mathcal {G}_n$ acts on
the left of $L^2(\mathbb {R}^n)$ according to the rule $(x,y) \,
f(t) = e^{2\pi iy \cdot t} f(t+ x)$, and extending to
the whole of $\mathbb {C}*\mathcal {G}_n$ by
$\mathbb {C}$-linearity.  To check that this indeed defines an action,
we need only verify that $(a,b)\bigl((x,y) f(t) \bigr) = \bigl(
(a,b) (x,y) \bigr) f(t) $, which is indeed true because both sides
equal $e^{2\pi i (a \cdot y  + b \cdot t + y \cdot t)} f(t+a+x)$.
Thus we obtain a homomorphism from $\mathbb {C}*\mathcal {G}_n$
into $\mathcal B(L^2(\mathbb {R}^n))$.
Since $\mathbb {C}*\mathcal {G}_n$ is a simple ring by
Lemma~\ref{Lsimple}, this homomorphism must be a monomorphism
and so
we may view $\mathbb {C}*\mathcal {G}_n$ as a $\mathbb
{C}$-subalgebra of $\mathcal B(L^2(\mathbb {R}^n))$.  We shall
consider the following conjecture.
\begin{Conj} \label{Cconj}
Let $0 \ne \theta \in \mathbb {C}*\mathcal {G}_n$ and
$0 \ne f\in L^2(\mathbb {R}^n)$.
Then $\theta f \ne 0$.
\end{Conj}

Motivation for studying this problem comes from Gabor analysis and
in particular the conjecture on page 2790 of \cite{heilram}.
If $G \leqslant \mathcal {G}_n$, then $\mathbb {C}*G$ will
denote the $\mathbb {C}$-subalgebra of
$\mathbb {C}*\mathcal {G}_n$ which has $\mathbb {C}$-basis $\{\bar g
\mid g \in G\}$.  Of course when talking about discrete subsets of
$\mathcal {G}_n$, we are giving $\mathcal {G}_n$ the usual topology
from $\mathbb {R}^{2n}$.
We shall prove
\begin{Thm}  \label{Tfirst}
Let $G$ be a discrete subgroup of $\mathcal {G}_n$.
If $0 \ne \theta \in \mathbb {C}*G$ and $0 \ne f\in
L^2(\mathbb {R}^n)$, then $\theta f \ne 0$.
\end{Thm}

Of course it follows immediately that if $G$ is a discrete subgroup of
$\mathcal {G}_n$, $g \in \mathcal {G}_n$,
$0 \ne \theta \in g\mathbb {C}*G$ and
$0 \ne f \in L^2(\mathbb {R}^n)$, then $\theta f \ne 0$.  This means
we can rephrase the above result in terminology closer to that of
\cite{heilram} as follows.
For $x,y \in \mathbb {R}^n$ and $f \in L^2(\mathbb {R}^n)$, define
$(x,y) f(t) = e^{2\pi iy \cdot t} f(t+x)$ and
if $\Lambda \subseteq \mathbb
{R}^{2n}$, define $S(f, \Lambda) = \{(x,y)f \mid (x,y) \in
\Lambda \}$.  Then Theorem~\ref{Tfirst} yields
\begin{Prop} \label{Pfirst}
Let $n$ be a positive integer, let
$\Lambda$ be a subset of  $\mathbb {R}^{2n}$ of the form $g + G$ where
$G$ is a discrete subgroup of $\mathbb {R}^{2n}$, and let
$0 \ne f \in L^2(\mathbb {R}^n)$.  Then $S(f,\Lambda)$ is linearly
independent.
\end{Prop}

As a byproduct, we shall obtain results on the von Neumann algebra
generated by $\mathbb {C}*G$, which we shall denote by $W*G$.
Thus $W*G$ is the
weak closure of $\mathbb {C}*G$ in $\mathcal B(L^2(\mathbb {R}^n))$
and is rather similar to the group von Neumann algebra of
$G$.  For $f,g \in L^2(\mathbb {R}^n)$, let
$\langle f,g \rangle$ denote the inner product $\int_{\mathbb {R}^n}
f(t) \bar g(t) \, dt$, where $\bar { }$ denotes complex conjugation,
and let $\mathcal U(L^2(\mathbb {R}^n))$ denote the set of closed
densely defined linear operators \cite[\S 2.7]{kadison} acting
on $L^2(\mathbb {R}^n)$.
Then the adjoint $\alpha^*$ of $\alpha \in \mathcal U(L^2(\mathbb
{R}^n))$ satisfies $\langle \alpha f, g \rangle = \langle f,
\alpha^* g \rangle$ whenever $f,g \in L^2(\mathbb {R}^n)$ and $\alpha
f$, $\alpha^* g$ are defined.  Of course $^*$ restricts to an
involution on both $\mathcal B(L^2(\mathbb {R}^n))$ and $W*G$.
If $G$ is a discrete subgroup of $\mathcal {G}_n$, then $W*G$ is a
\emph{finite} von Neumann algebra by Lemma~\ref{Ltrace}; also in many
cases this can be deduced from Rieffel's paper \cite{rieffel}.
In this situation, we let
$U*G$ indicate the operators of $\mathcal U(L^2(\mathbb {R}^n))$
which are affiliated to $W*G$ \cite[p.~150]{ber}.
The results of \cite{ber} (especially theorem~1 and
the proof of theorem~10) now show that $(U*G)^* = U*G$, $U*G$ is a
$*$-regular ring containing $W*G$, and every element of $U*G$ can be
written in the form $\gamma \delta^{-1}$ where $\gamma, \delta \in
W*G$.  In particular every nonzero divisor in $W*G$ is
invertible in $U*G$.  We shall prove
\begin{Thm} \label{Tmain}
Let $G$ be a discrete subgroup of $\mathcal {G}_n$.  Then
$W*G$ is a finite von Neumann algebra,
every nonzero element of $\mathbb {C}*G$ is invertible in $U*G$, and
the set $\{\gamma \delta^{-1} \mid \gamma \in \mathbb {C}*G,
0 \ne \delta \in \mathbb {C}*G \}$ is a division subring of $U*G$.
\end{Thm}

Let $L$ be a locally compact group,
let $G$ be a torsion free subgroup of $L$, and let $L^2(L)$ denote the
Hilbert space of square integrable functions on $L$ with respect to
the left Haar measure on $L$.  Then $G$ acts on the left of $L^2(L)$
according to the rule $g f(l) = f(g^{-1}l)$ for $g\in G, f \in L^2(L),
l \in L$.  For $f \in L^2(L) \setminus 0$, a closely related problem
to Conjecture~\ref{Cconj} is to determine whether
the set $\{gf \mid g\in G \}$ is
linearly independent over $\mathbb {C}$.  If
the von Neumann algebra $W*G$
generated by $G$ is a finite von Neumann algebra, then by using the
techniques of this paper, it is possible in many cases to show that
the set $\{gf \mid g \in G \}$ is linearly independent.  On the other
hand if  $W*G$ is not a finite von Neumann algebra, then the
techniques of this paper cannot be applied.  It will usually be the
case that $W*G$ is not finite if $G$ is not
discrete and has no abelian subgroup of finite index.  A specific
example would be to let $L$ be the Heisenberg group consisting of
upper unitriangular 3 by 3 matrices with entries in $\mathbb {R}$,
in other words matrices of the form
\[
\begin{pmatrix}
1&a&b\\
0&1&c\\
0&0&1
\end{pmatrix}
\]
where $a,b,c \in \mathbb {R}$,
and to let $G = L$.
Then it is not known in this case whether for $f \in L^2(L) \setminus
0$, the set $\{gf \mid g \in G\}$ is linearly independent.

I am very grateful to Chris Heil and Jonathan Rosenblatt for bringing
the problem studied in this paper to my attention, and for some
helpful email correspondence.

\section{Notation, Terminology and Assumed Results}

The identity of a group will be denoted by either 0 or 1.
If $n$ is a positive integer and $R$ is a ring, then $\Mat_n(R)$ will
denote the $n$ by $n$ matrices over $R$,
and we shall let $\delta_{ij}$
indicate the Kronecker delta, so $\delta_{ij} = 0$ if $i\ne j$ and
$\delta_{ij} = 1$ if $i=j$.
The identity matrix of $\Mat_n(R)$ will be denoted by $I_n$,
and the zero matrix of $\Mat_n(R)$ will be denoted by $0_n$.
We shall view vectors in $\mathbb {R}^n$
as column vectors rather than row vectors.
A lattice in $\mathbb {R}^n$ will mean a
discrete subgroup of $\mathbb {Z}$-rank $n$; in other words a discrete
subgroup of finite covolume (note that this is a different definition
of lattice from that of \cite[p.~2791]{heilram}).  If $\alpha =
\sum_{g\in \mathcal {G}_n} \lambda_g g \in \mathbb {C}*\mathcal
{G}_n$ where $\lambda_g \in \mathbb {C}$
for all $g\in \mathcal {G}_n$, then the support of $\alpha$,
denoted $\supp \alpha$, is the  set
$\{g\in \mathcal {G}_n \mid \lambda _g \ne 0\}$.
We shall use the notation $\| f\|_2$ for the norm $\sqrt{\langle f,
f \rangle}$ of an element $f \in L^2(\mathbb {R}^n)$,
and $\overline{X}$
for the closure of a subset $X$ in $L^2(\mathbb {R}^n)$.
The commutant of
a subset $A$ of $\mathcal B(L^2(\mathbb {R}^n))$ is $A' = \{x
\in \mathcal B(L^2(\mathbb {R}^n)) \mid ax = xa$ for all $a \in A \}$.
If $A = A^*$, then  $A'$ is a von Neumann algebra and
by von Neumann's double commutant
theorem \cite[theorem~1.2.1]{arveson}, $A$ is dense in $A''$ in the
weak operator topology.
Thus another description of $W*G$ is the double commutant
of $\mathbb {C}*G$ in $\mathcal B(L^2(\mathbb {R}^n))$.
In the case $W*G$ is a finite von Neumann algebra, we can now
describe $U*G$ as those unbounded operators in $\mathcal U(L^2(\mathbb
{R}^n))$ which commute with every element of $(W*G)'$.

\begin{Lem} \label{Lsimple}
$\mathbb {C}*\mathcal {G}_n$ is a simple ring.
\end{Lem}
\begin{proof}
Suppose $0 \ne I \lhd \mathbb {C}*\mathcal {G}_n$ with
$I \ne \mathbb {C}*\mathcal {G}_n$,
and choose $0 \ne \alpha \in I$ with
minimal support.  If $g \in \supp \alpha$, then $1 \in \supp \bar
g^{-1} \alpha$ and $\bar g^{-1}\alpha \in I$, so we may assume that
$1\in \supp \alpha$.  Since $I \ne \mathbb {C}*\mathcal {G}_n$, we
may choose $a\in \mathcal {G}_n$ such
that $1 \ne a \in \supp \alpha$.  Then there
exists $g \in \mathcal {G}_n$ such that
$\bar g \bar a \bar g^{-1} \ne \bar a$,
and now we have $0 \ne \bar g \alpha \bar g^{-1} - \alpha \in I$.
This  contradicts the minimality of $\supp \alpha $ because
$|\supp (\bar g\alpha \bar g^{-1} -\alpha)| < |\supp \alpha |$,
and the result follows.
\end{proof}

If $R$ is a ring and $\sigma $ is an automorphism of $R$, then
$R_\sigma [X]$ will denote the twisted polynomial ring over $R$ in the
indeterminate $X$, so multiplication is defined by $\sum a_iX^i \sum
b_j X^j = \sum_n (\sum_{i+j=n} a_i \sigma^i b_j )X^n$.
We say that $R$ is an Ore domain if it is contained in a division ring
$D$, called the division ring of fractions of $R$,
such that every element of $D$ can be written in the form $rs^{-1}$
and also in the form $s^{-1}r$, with $r,s \in R$ and $s\ne 0$.
Of course the division ring $D$ containing $R$ is unique up to
$R$-isomorphism.  Also if $R$ is contained in a ring $D'$ such that
every nonzero element of $R$ is invertible, then the set $\{rs^{-1}
\mid r,s \in R$ and $s \ne 0\}$ is the division ring of fractions
containing $R$.  The following two elementary results are well known.
\begin{Lem}  \label{Lore1}
Let $R$ be an Ore domain with division ring of fractions $D$, and let
$\sigma$ be an automorphism of $R$.  Then $\sigma$ extends uniquely to
an automorphism of $D$, which we shall also call $\sigma$, and if
$\alpha, \beta \in D_{\sigma}[X]$, then there exists
$r \in R\setminus 0$ such that $r\alpha, r\beta \in R_{\sigma}[X]$.
\end{Lem}

\begin{Lem} \label{Lore2}
Let $G$ be a subgroup of $\mathcal {G}_n$.
Then $\mathbb {C}*G$ is an Ore
domain, and if $I,J$ are nonzero left ideals of $\mathbb {C}*G$, then
$I\cap J \ne 0$.
\end{Lem}

Finally we require the following:

\begin{Lem} \label{Lbilinear}
Let $G$ be a discrete subgroup of $\mathcal {G}_n$, let $H \lhd G$
such that $G/H$ is infinite cyclic, and let $x \in G$ such that $Hx$
is a generator for $G/H$.
If $\zeta \in \mathbb {C}$ and $|\zeta| =1$,
then there exists $y\in \mathcal {G}_n$ such that $\bar y \bar h \bar
y^{-1} =\bar h $ for all $h \in H$ and $\bar y \bar x \bar y^{-1} =
\zeta \bar x$ in $\mathbb {C}*\mathcal {G}_n$.
\end{Lem}
\begin{proof}
Since $G$ is discrete, we may choose $m\in \mathbb {Z}$ and a subset
$\{h_1, \dots, h_m\}$  which generates $H$ and is linearly independent
over $\mathbb {R}$.
Note that $\{h_1, \dots, h_m, x\}$ is also linearly
independent over $\mathbb {R}$.  Choose $t\in \mathbb {R}$ such that
$e^{2\pi i t} = \zeta$, and define
a bilinear form $\beta \colon \mathcal {G}_n \to \mathbb {R}$ by
$\beta \bigl( (a,b), (c,d) \bigr)
= a \cdot d - b \cdot c$, where $a,b,c,d \in \mathbb {R}^n$.
Note that in $\mathbb
{C}*\mathcal {G}_n$, we have
\[
(a,b)(c,d)(a,b)^{-1} = e^{2\pi i (a \cdot d -b \cdot c)} (c,d).
\]
It is easily checked that $\beta$ is
nondegenerate, so there exists $y \in \mathcal {G}_n$ such that $\beta
(y,h_i) = 0$ for all $i$ and $\beta (y,x) = t$.
This completes the proof.
\end{proof}

\section{Faithful Traces}

In this section, we show that $W*G$ has a faithful weakly continuous
tracial state, which in particular will establish that $W*G$ is a
finite von Neumann algebra.
Throughout this section, $n$ will be a positive
integer.  The purpose of the next lemma is to reduce to the case when
$G$ is a lattice in $\mathbb {R}^{2n}$ containing $1 \times \mathbb
{Z}^n$; its proof is modelled on \cite[\S 2, p.~2790]{heilram}.

We shall think of $\mathbb {R}^{2n}$ as $\mathbb {R}^n \oplus
\mathbb {R}^n$, so we can view $\mathbb {R}^n$ as a subgroup of
$\mathbb {R}^{2n}$ in the usual way via the map $x \mapsto (x,0)$.
We then have a monomorphism
$\psi \colon \mathcal {G}_n \to \mathcal {G}_{2n}$ and this
induces a monomorphism $\mathbb {C} * \mathcal {G}_n \to \mathbb {C} *
\mathcal {G}_{2n}$, which we shall also call $\psi$.

Given $f,g \in L^2(\mathbb {R}^n)$, we can form the element $f\otimes
g\in L^2(\mathbb {R}^{2n})$
defined by $(f \otimes g) (x,y) = f(x)g(y)$
for $x,y \in L^2(\mathbb {R}^n)$, and then the functions of the form
$\sum_{i=1}^m f_i \otimes g_i$
are dense in $L^2(\mathbb {R}^{2n})$.  If $\theta \in
\mathcal {B} (L^2(\mathbb {R}^n))$, then we have a well defined
operator $\theta \otimes 1 \in \mathcal {B} (L^2(\mathbb {R}^{2n}))$
satisfying $(\theta \otimes 1) (f \otimes g) = (\theta f) \otimes g$
for all $f,g \in L^2(\mathbb {R}^n)$, and this yields a weakly
continuous $*$-monomorphism $\theta \mapsto \theta \otimes 1  \colon
\mathcal {B}(\mathbb {R}^n) \to \mathcal {B}( \mathbb {R}^{2n})$.

Note that when we view $\mathbb {C}*\mathcal {G}_n$ and $\mathbb {C}*
\mathcal {G}_{2n}$ as subalgebras of $\mathcal {B} (L^2(\mathbb
{R}^n))$ and $\mathcal {B} (L^2(\mathbb {R}^{2n}))$
respectively, then $\psi
(\theta ) = \theta \otimes 1 $ for all $\theta \in \mathbb {C} *
\mathcal {G}_n$.  Furthermore if $G \le \mathcal {G}_n$, then
$\psi$ induces isomorphisms $W*G \to W* \psi G$ and
(assuming $W*G$ is a finite von Neumann algebra) $U*G \to U * \psi
G$, which means we may identify $G$ with the subgroup $\psi G$ of
$\mathcal {G}_{2n}$; we shall do this without further comment and
without using $\psi$ in the future.

Let $\{e_1, \dots, e_{2n} \}$ denote the standard basis for $\mathbb
{R}^{2n}$, so $e_i$ has a 1 in the $i$th position and zeros elsewhere,
and $e_i \cdot e_j = \delta _{ij}$.
If $G \leqslant \mathcal {G}_n$, then we define
$\{\mathbb {C}G\} = \{ \lambda g \mid \lambda \in \mathbb {C}\text{
and } g \in G\}$, a subset of $\mathbb {C}*G$.

\begin{Lem} \label{Lmetaplectic}
Let $G$ be a discrete subgroup of $\mathcal {G}_n$.  Then there exists
a lattice $H$ in $\mathcal {G}_{2n}$ containing $1 \times \mathbb
{Z}^{2n}$ and a unitary operator $u \in \mathcal {B}(L^2(\mathbb
{R}^{2n}))$, such that $u \{ \mathbb {C}G \} u^{-1} \subseteq \{
\mathbb {C}H\}$.
\end{Lem}
\begin{proof}
Choose an $\mathbb {R}$-basis $\{g_1, \dots, g_{2n}\}$ for $\mathcal
{G}_n$
such that $\{g_1, \dots , g_r\}$ is a $\mathbb {Z}$-basis for $G$,
where $r$ is the rank of $G$.
Let $\mathcal {E} = \{ (e_1,0), \dots, (e_{2n},0), (0,e_1), \dots,
(0,e_{2n}) \}$, let
\begin{multline*}
\mathcal {F} = \{
(e_1,e_{n+1})/\sqrt 2, (e_2,e_{n+2})/\sqrt 2, \dots,
(e_n,e_{2n})/\sqrt 2,
(e_{n+1},e_1)/\sqrt 2,\\ (e_{n+2}, e_2)/\sqrt 2, \dots,
(e_{2n},e_n)/\sqrt 2,
(-e_{n+1},e_1)/\sqrt 2,
(-e_{n+2}, e_2)/\sqrt 2, \dots,\\
(-e_{2n},e_n)/\sqrt 2,
(-e_1,e_{n+1})/\sqrt 2, (-e_2,e_{n+2})/\sqrt 2, \dots,
(-e_n,e_{2n})/\sqrt 2 \},
\end{multline*}
and let
\begin{multline*}
\mathcal {K} = \{
g_1, \dots, g_{2n},
(-e_{n+1},e_1)/\sqrt 2,
(-e_{n+2}, e_2)/\sqrt 2, \dots,
(-e_{2n},e_n)/\sqrt 2,\\
(-e_1,e_{n+1})/\sqrt 2, (-e_2,e_{n+2})/\sqrt 2, \dots,
(-e_n,e_{2n})/\sqrt 2 \},
\end{multline*}
so $\mathcal{E}$, $\mathcal {F}$ and $\mathcal {K}$ are
$\mathbb {R}$-bases of $\mathcal {G}_{2n}$.
For $i=1, \dots, 4n$, we shall let $\hat e_i$, $f_i$, $k_i$ denote
the $i$th basis elements of $\mathcal {E}$, $\mathcal {F}$, $\mathcal
{K}$ respectively, and we shall let $K$ be the lattice in $\mathcal
{G}_{2n}$ which has $\mathbb {Z}$-basis $\mathcal {K}$.
Let $A_i$ denote the coordinates of $k_i$
with respect to the basis $\mathcal {F}$, and let $a_{ji}$ denote the
$j$th coordinate of $A_i$.  Then for $2n + 1 \le i \le 4n$,
$a_{ji}=1$ if $j=i$ and $a_{ji}=0$ if $j\ne i$.
Now define
$h_i = \sum_{j=1}^{4n} a_{ji} \hat e_j \in \mathcal {G}_{2n}$,
and let $H$ be the subgroup of $\mathcal
{G}_{2n}$ generated by the $h_i$.  Then $H$ is a lattice in $\mathcal
{G}_{2n}$ which contains $1 \times \mathbb {Z}^{2n}$.

Let $T$ be the
transition matrix from $\mathcal {E}$ to $\mathcal {F}$, and let $J =
J_{2n} =
\begin{pmatrix}
0_n & I_n\\
I_n & 0_n
\end{pmatrix}
\in \Mat_{2n}(\mathbb {R})$.  Thus if $T$
has entries $t_{ij}$, then $f_j = \sum_{i=1}^{4n} t_{ij} \hat e_i$,
and if we
think of the $A_i$ as column vectors, then the coordinates of
$k_i$ with respect to $\mathcal {E}$ are $T A_i$.  Also
\begin{align*}
T &=
\begin{pmatrix}
I_{2n}/\sqrt 2 & -J_{2n}/\sqrt 2 \\
J_{2n}/\sqrt 2 & I_{2n}/\sqrt 2
\end{pmatrix}\\
&=\begin{pmatrix}
I_{2n} & 0_{2n}\\
-J_{2n} & I_{2n}
\end{pmatrix}
\begin{pmatrix}
J_{2n}/\sqrt 2 & 0_{2n} \\
0_{2n} & J_{2n}\sqrt 2
\end{pmatrix}
\begin{pmatrix}
0_{2n}&-I_{2n}\\
I_{2n} & 0_{2n}
\end{pmatrix}
\begin{pmatrix}
I_{2n} & 0_{2n}\\
-J_{2n} & I_{2n}
\end{pmatrix}.
\end{align*}

Let $\tau, \alpha, \beta, \gamma
\colon \mathcal {G}_{2n}\to \mathcal {G}_{2n}$ be the linear
mappings determined by the matrices
\[
T,\
\begin{pmatrix}
J_{2n}/\sqrt 2 & 0_{2n} \\
0_{2n} & J_{2n}\sqrt 2
\end{pmatrix},\
\begin{pmatrix}
I_{2n} & 0_{2n}\\
-J_{2n} & I_{2n}
\end{pmatrix},
\begin{pmatrix}
0_{2n} & -I_{2n}\\
I_{2n} &0_{2n}
\end{pmatrix}
\]
respectively with respect to the basis $\mathcal {E}$, so $\tau \hat
e_i = f_i$ for all $i$.
Then $\tau H = K \supseteq G$, so it will be sufficient
to show that there exists a unitary operator $u
\in \mathcal {B}(L^2(\mathbb {R}^{2n}))$ such that $u^{-1} \mathbb {C}
g u = \mathbb {C} \tau g$ for all $g \in \mathcal {G}_{2n}$.
Since $\tau = \beta \alpha
\gamma \beta$, it will be sufficient to do this with $\alpha, \beta,
\gamma$ in  place of $\tau$.  We now use metaplectic transformations
\cite[p.~578]{stein}.
Write $g = (x,y)$ where $x,y \in \mathbb
{R}^{2n}$, and then we have three cases to consider.
\begin{enumerate}
\item
The matrix $\alpha =
\begin{pmatrix}
J_{2n}/\sqrt 2 & 0_{2n}\\
0_{2n} & \sqrt 2 J_{2n}
\end{pmatrix}$.  For $f \in L^2 (\mathbb {R}^{2n})$ and $t \in \mathbb
{R}^{2n}$, we define $uf(t) = 2^{-n/2}
f(Jt/\sqrt 2)$ (we are considering $t$ as a column vector in $\mathbb
{R}^{2n}$ here).  Then $u$ is $\mathbb
{C}$-linear and $\| uf\|_2 = \|f\|_2$ for all $f \in L^2(\mathbb
{R}^{2n})$, hence $u$ is a unitary operator.  Also $u^{-1}f(t) =
2^{n/2}f(\sqrt 2 Jt)$ because $J_{2n}^2 = I_{2n}$,
consequently
\begin{align*}
u^{-1} g uf(t) &= u^{-1} g 2^{-n/2} f(Jt/\sqrt 2)
= u^{-1} e^{2\pi i y \cdot t}2^{-n/2} f(J(t+x) /\sqrt 2) \\
&= e^{2 \pi i \sqrt 2 Jy \cdot t} f(t + Jx/\sqrt 2)
\quad \text{ because } J_{2n} \text{ is symmetric}\\
&= (\alpha g)f(t)
\end{align*}
for all $t \in \mathbb {R}^{2n}$ and for all $f \in
L^2 (\mathbb {R}^{2n})$.  Thus $u^{-1}gu = \alpha g$ as required.

\item The matrix $\beta =
\begin{pmatrix}
I_{2n} & 0_{2n}  \\
-J_{2n} & I_{2n}
\end{pmatrix}$.
Here we define $uf(t) = e^{-\pi iJt \cdot t} f(t)$.
Then $u$ is $\mathbb
{C}$-linear and $\|u f\|_2 = \|f \|_2$, so $u$ is a unitary operator.
Since $u^{-1} f(t) = e^{\pi iJt \cdot t}f(t)$,
\begin{align*}
u^{-1} g uf(t) &= u^{-1}
g e^{-\pi iJt \cdot t} f(t) = u^{-1} e^{2\pi i y \cdot t}e^{-\pi
iJ(t+x) \cdot (t+x)} f(t+x)\\
&= e^{\pi iJt \cdot t} e^{2\pi i y \cdot t} e^{-\pi
iJ(t+x) \cdot (t+x)} f(t+x)\\
& = e^{-\pi i Jx \cdot x} e^{2\pi it \cdot (y-Jx)} f(t+x) \quad \text{
because } J_{2n} \text{ is symmetric}\\
&=e^{-\pi i Jx \cdot x} (\beta g) f(t)
\end{align*}
for all $t\in \mathbb {R}^{2n}$, and we have shown that $u^{-1} gu \in
\mathbb {C} (\beta g)$.

\item The matrix $\gamma =
\begin{pmatrix}
0_{2n}&-I_{2n}\\
I_{2n}&0_{2n}
\end{pmatrix}$.  Here we use the Fourier transform; specifically
$uf(t) = \int_{\mathbb {R}^{2n}} e^{2\pi i t \cdot s} f(s)\,ds$,
and then $u^{-1} f(t) = \int_{\mathbb {R}^{2n}} e^{-2\pi it \cdot s}
f(s)\,ds$.
Observe that
$(x,0) u f(t) = u(0,x)f(t)$ and $(0,y)uf(t) = u (-y,0)f(t)$,
consequently
\[
u^{-1}(x,y)u = u^{-1}(0,y)(x,0)u = (-y,0) (0,x) = e^{-2\pi i
x \cdot y} (-y,x) = e^{-2\pi i x \cdot y}\gamma (x,y)
\]
and we deduce that $u^{-1} \mathbb {C} g u = \mathbb {C} \gamma g$.
\end{enumerate}
This completes the proof of Lemma~\ref{Lmetaplectic}
\end{proof}

\begin{Lem}  \label{Ltrace}
Let $G$ be a discrete subgroup of $\mathcal {G}_n$, and
define $\tau \colon \mathbb {C}*G \to \mathbb
C$ by $\tau g = 0$ when $1 \ne g \in G$, and $\tau 1 = 1$.  Then
\begin{enumerate}[\normalfont (i)]
\item $\tau$ extends to a weakly continuous $\mathbb {C}$-linear
map $W*G \to \mathbb {C}$.
\label{Ltrace1}

\item If $\alpha, \beta \in W*G$, then $\tau (\alpha \beta) =
\tau (\beta \alpha ) $.
\label{Ltrace2}

\item If $\alpha \in W*G$ and $x \in \mathcal G_n$, then
$\tau (\bar x \alpha \bar x^{-1} ) = \tau (\alpha)$.  \label{Ltrace21}

\item If $e$ is a nonzero projection in $W*G$, then $0 < \tau e \le
1$.
\label{Ltrace3}

\item Let $e,f$ be projections in $W*G$.  If $e L^2(\mathbb {R}^n)
\subseteq fL^2(\mathbb {R}^n)$, then $\tau e \le \tau f$.
\label{Ltrace4}

\item Let $e,f$ be projections in $W*G$, and let $h$ be the projection
of $L^2(\mathbb {R}^n)$ onto $\overline{e L^2(\mathbb {R}^n)
+ f L^2(\mathbb
{R}^n)}$.  Then $h \in W*G$ and if $eL^2(\mathbb {R}^n) \cap
fL^2(\mathbb {R}^n) = 0$, then $\tau e + \tau f = \tau h$.
\label{Ltrace5}
\end{enumerate}
\end{Lem}
\begin{proof}
Since $G$ is a discrete subgroup of $\mathcal {G}_{n}$,  there is
by Lemma~\ref{Lmetaplectic}
a lattice $H$ in $\mathcal {G}_{2n}$ containing $1 \times \mathbb
{Z}^{2n}$ and a unitary operator $u \in
\mathcal B(L^2(\mathbb {R}^{2n}))$ such that
$u \mathbb {C}*G u^{-1} \subseteq \mathbb {C}*H$.  If
we can find a weakly continuous $\mathbb {C}$-linear map $\tau \colon
\mathbb {C}*H \to \mathbb {C}$ with the required properties, then
the weakly continuous $\mathbb {C}$-linear map $\alpha \mapsto \tau
(u\alpha u^{-1})$ for $\alpha \in \mathbb {C}*G$
will suffice.  Therefore we may assume that $G$
is a lattice in $\mathcal {G}_n$ containing $1 \times \mathbb {Z}^n$.
If $a$ is a positive number, we shall let $\mathcal C(a)$ denote the
standard unit cube in $\mathbb {R}^n$ with side of length $a$; thus
$\mathcal C(a) = \{(a_1, \dots, a_n ) \mid 0 \le a_i \le a$ for all
$i \}$.

\eqref{Ltrace1}
Choose a positive integer $b$ such that $h \mathcal C(1/b)
\cap \mathcal C(1/b) = \emptyset$ whenever $(h,k) \in G \setminus
(1 \times \mathbb {Z}^n)$, which is possible because $G$ is a lattice
in $\mathcal {G}_n$ containing $1 \times \mathbb {Z}^n$,
and set $\mathcal C = \mathcal C(1/b)$.  Let $c = b^n$, let
$\mathcal C_1, \dots,
\mathcal C_c$ denote the $c$ translates of $\mathcal  C$ which
are contained in the unit cube $\mathcal C(1)$, and for each $i$, let
$\chi_i $ denote the characteristic function of $\mathcal C_i$.
For $\theta \in W*G$, define
\[
\tau \theta = \sum_{i=1}^c \langle \theta \chi_i, \chi_i \rangle =
\sum_{i=1}^c \int_{\mathcal C_i} \theta \chi_i (t) \, dt.
\]
Let $g \in G$ and write $g = (h,k)$ where $h,k \in \mathbb
R^n$.
Then $g \chi_i(t) = e^{2\pi i k \cdot t}\chi_i(t + h)$, so if $h \ne
0$ we have $g \chi_i(t) = 0$ for all $t \in \mathcal C_i$ and hence
$\tau g =0$.  On the other hand if $h = 0$ and $k \ne 0$, then
$k \in \mathbb {Z}^n \setminus 0$ because $G$ is a lattice
containing $1 \times \mathbb {Z}^n$, consequently
\[
\tau g = \sum_{i=1}^c \int_{\mathcal C_i} e^{2\pi i k \cdot t} \, dt
= \int_{\mathcal C} e^{2 \pi i k \cdot t} \, dt = 0.
\]
Finally $\tau 1 = 1$ and \eqref{Ltrace1} is proven.

\eqref{Ltrace2}
If $x,y \in G$, then $xy = 1$ if and only if $yx = 1$.  Therefore
$\tau \bar x\bar y = \tau \bar y\bar x = 0$ if $xy \ne 1$ and
$\bar x\bar y = \bar y\bar x$ if $xy = 1$,
hence $\tau \bar x\bar y = \tau \bar y\bar x$ for all $x,y
\in G$ and we deduce that $\tau \alpha\beta = \tau \beta \alpha$ for
all $\alpha, \beta \in \mathbb {C}*G$.  Since $\tau$ is weakly
continuous and $W*G$ is the weak closure
of $\mathbb {C}*G$, we see that
$\tau \alpha \beta = \tau \beta \alpha$ for all $\alpha, \beta \in
W*G$ which proves \eqref{Ltrace2}.

\eqref{Ltrace21}
Define $\sigma \colon W*G \to \mathbb {C}$ by $\sigma (\alpha) =
\tau (x \alpha x^{-1})$.  Observe that $\sigma 1 = 1$
and if $1 \ne g \in G$, then
$\sigma g = \tau (xgx^{-1}) = 0$ because $xgx^{-1} = \zeta g$ for some
$\zeta \in \mathbb {C}$ with $| \zeta | = 1$.
Thus $\sigma g = \tau g$
for all $g \in G$ and since $\sigma$ is a weakly continuous $\mathbb
C$-linear map, we deduce that $\sigma (\alpha ) = \tau (\alpha) $ for
all $\alpha \in W*G$, which is the required result.

\eqref{Ltrace3}
Note that if $e$ is a projection in $W*G$, then $e^*e = e$ and hence
\[
\tau e = \sum_{i=1}^c \langle e^* e \chi_i, \chi_i \rangle =
\sum_{i=1}^c \langle e\chi_i, e \chi_i \rangle \ge 0.
\]
Let $K =
\mathbb {Z}^n \times \mathbb {Z}^n \leqslant \mathcal G_n$.  If
$k\in K$, then $\tau e = \tau (k^{-1} e k)$ by \eqref{Ltrace21} and we
deduce that
\[
\tau e = \tau (k^{-1} ek) = \sum_{i=1}^c \langle k^{-1} e^*e k \chi_i,
\chi_i \rangle = \sum_{i=1}^c \langle ek \chi_i, ek \chi_i \rangle.
\]
Let $\chi$ denote the characteristic function of
$\mathcal C(1)$, and suppose $\tau e = 0$.  Then $ ek \chi_i = 0$ for
all $i$, hence $ek \chi = 0$ for all $k \in K$.  Now the set $\{ k\chi
\mid k \in K \}$ forms a Hilbert basis for $L^2(\mathbb {R}^n)$ so if
$\tau e = 0$, we see that $ef = 0$ for all $f \in L^2(\mathbb {R}^n)$
and we deduce that $e=0$.  Also $1-e$ is a projection if $e$ is a
projection, so applying the above to $1-e$ we obtain $0 \le \tau
(1-e)$, hence $\tau e \le 1$ and \eqref{Ltrace3} follows.

\eqref{Ltrace4}
Let $h$ be the projection of $L^2(\mathbb {R}^n)$ onto the orthogonal
complement of $eL^2(\mathbb {R}^n)$ in $fL^2(\mathbb {R}^n)$.  Then
$e + h = f$, hence $\tau e + \tau h = \tau f$.  Thus $h\in W*G$ and
the result follows from \eqref{Ltrace3}.

\eqref{Ltrace5}
Let $u$ be a unitary operator in $(\mathbb {C}*G)'$.  Then $ueu^{-1}
= e$ and $ufu^{-1} = f$.  Since $h$ is the projection of $L^2(\mathbb
R^n)$ onto $\overline{eL^2(\mathbb {R}^n) + f L^2(\mathbb {R}^n)}$, we
see that $uhu^{-1}$ is the projection of $L^2(\mathbb {R}^n)$ onto
$\overline{ueu^{-1} L^2(\mathbb {R}^n) +
ufu^{-1} L^2(\mathbb {R}^n)} =
\overline{eL^2(\mathbb {R}^n) + fL^2(\mathbb {R}^n)}$ and we deduce
that $uhu^{-1} = h$.  Therefore $uh = hu$.
Now $(\mathbb {C}*G)'$ is a von Neumann algebra, so any
element of $(\mathbb {C}*G)'$
is a $\mathbb {C}$-linear sum of unitary elements, hence
$xh = hx$ for all $x \in (\mathbb {C}*G)'$ and we conclude that $h
\in W*G$.

We now claim that $h = e \cup f$ \cite[p.~4]{berberian}.  Since
$eL^2(\mathbb {R}^n) \subseteq hL^2(\mathbb {R}^n)$,
we see that $e = he$
and hence $e \leq h$.  Similarly $f \le h$ and so $e \cup f \le h$.
Now let $g = e \cup f$.  Then $gL^2(\mathbb {R}^n)
\supseteq eL^2(\mathbb
R^n), fL^2(\mathbb {R}^n)$ and hence $gL^2(\mathbb {R}^n) \supseteq
hL^2(\mathbb {R}^n)$.  We deduce that $g \ge h$, consequently $g=h$
and the claim is established.

If $eL^2(\mathbb {R}^n) \cap fL^2(\mathbb {R}^n) = 0$, then $(e \cap
f)L^2(\mathbb {R}^n) \subseteq eL^2(\mathbb {R}^n)
\cap fL^2(\mathbb {R}^n)
= 0$.  Since $W*G$ is a von Neumann algebra, we may apply the
parallelogram law to deduce that $e \sim e\cup f - f$ \cite[\S 1,\S
13]{berberian}.  Thus there is an element $w \in W*G$ such that $w^*w
= e$ and $ww^* = e \cup f - f$.  Since $\tau (w^*w) = \tau (w w^*)$,
we deduce that $\tau e = \tau(e \cup f - f)$ and hence $\tau e + \tau
f = \tau (e \cup f) = \tau h$.  This completes the proof.
\end{proof}

\section{Proofs}
Theorems~\ref{Tfirst} and \ref{Tmain} are now immediate consequences
of the following result.

\begin{Lem}  \label{Lmain}
Let $G$ be a discrete subgroup of $\mathbb {R}^n$ and let $\theta \in
\mathbb {C}*G \setminus 0$.  Then
\begin{enumerate}[\normalfont (i)]
\item If $0 \ne f \in L^2(\mathbb {R}^n)$, then $\theta f \ne 0$.
\label{Lmain1}

\item $\theta$ is invertible in $U*G$.
\label{Lmain2}

\item The set $\{\gamma\delta^{-1} \mid
\gamma \in \mathbb {C}*G, 0 \ne
\delta \in \mathbb {C}*G\}$ is a division
subring of $U*G$, and is equal
to $\{ \delta^{-1} \gamma \mid \gamma \in \mathbb {C}*G,
0 \ne \delta \in \mathbb {C}*G \}$.
\label{Lmain3}
\end{enumerate}
\end{Lem}
\begin{proof}
Since $G$ is a discrete subgroup of $\mathcal {G}_n$,
it is a free abelian
group of rank at most $2n$.  We shall prove the result by induction on
the rank of $G$, the result being trivially true if the rank of $G$ is
zero, because then $G = 1$.  Thus we may assume that the rank of $G$
is strictly positive, and then there exists $H \lhd G$ such that
$G/H \cong \mathbb {Z}$.  Since $H$ has
strictly smaller rank than $G$,
we may assume that the result is true for $H$.  Let $\tau \colon W *G
\to \mathbb {C}$ be the weakly continuous tracial state obtained from
Lemma~\ref{Ltrace}.

\eqref{Lmain1}
For $\alpha \in \mathbb {C}*G$, let $\ker\alpha = \{f \in L^2(\mathbb
R^n) \mid \alpha f = 0\}$, and let $\mathcal N(\alpha)$
be the projection from
$L^2(\mathbb {R}^n) $ onto $\ker \alpha$.  Suppose $u$ is a unitary
element in $(\mathbb {C}*G)'$.  Then $u^{-1} \mathcal N(\alpha) u =
\mathcal N(u^{-1} \alpha u) = \mathcal N(\alpha)$.  Since
$(\mathbb {C}*G)'$  is a von Neumann algebra, every
element of $(\mathbb
C*G)'$ is a linear combination of unitary  elements of $(\mathbb
C*G)'$ and we deduce that $\mathcal N(\alpha) $ commutes with every
element of $(\mathbb {C}*G)'$.  Therefore $\mathcal N(\alpha ) \in
W * G$.

Let $\nu = \sup \{ \tau (\mathcal N(\alpha)) \mid 0 \ne
\alpha \in \mathbb {C}*G
\}$.  If $\nu = 0$ then $\mathcal N(\theta) = 0$
by Lemma~\ref{Ltrace}\eqref{Ltrace3},
hence $\ker \theta = 0$
and the result follows, so we may assume that $0 < \nu
\le 1$.  Therefore we may choose $\alpha \in \mathbb {C}*G$ such that
$\tau \mathcal N(\alpha ) > \nu /2$.  Since $G/H$ is infinite
cyclic, there exists $x \in G$ such
that $Hx$ generates $G/H$, and then
we may write $\alpha = \sum_{i=-\infty}^{i=\infty} \alpha_i x^i$ where
$\alpha_i \in \mathbb {C}*H$ and $\alpha_i = 0$ for all but finitely
many $i$.  By replacing $\alpha$ with $x^m \alpha$ for some
integer $m$, we may assume that $\alpha_1 \ne 0$ and $\alpha_i = 0$
for all $i< 0$.

By induction, there is a division subring $D$ of $U*H$ containing
$\mathbb {C}*H$ which is the division ring of fractions of
$\mathbb {C}*H$.  Let $\sigma$ be the automorphism $\beta
\mapsto x\beta x^{-1} \colon \mathbb {C}*H \to \mathbb {C}*H$.
By Lemma~\ref{Lore1} we may extend $\sigma$ to an automorphism of $D$,
which we shall also call $\sigma$.  We now have a natural ring
homomorphism $\theta \colon D_\sigma [X] \to U*G$,
defined by $\theta X
= x$ and $\theta d = d$ for all $d\in D$, which maps $(\mathbb
{C}*H)_\sigma [X] $ into $\mathbb {C} *G$.
By \cite[lemma~16]{zero}, there exists $\zeta \in \mathbb {C}$ with
$|\zeta| = 1$ and $\beta', \gamma' \in D_\sigma [X]$ such that
\[
\beta' \sum_i \alpha_0^{-1}\alpha_i X^i +
\gamma' \sum_i \alpha_0^{-1} \alpha_i \zeta^i X^i = 1.
\]
By Lemma~\ref{Lore1}, there exists $0 \ne r \in \mathbb {C}*H$
such that $r \beta' \alpha _0^{-1}, r \gamma' \alpha_0^{-1}
\in (\mathbb {C}*H)_\sigma [X]$, so setting $\beta = r \beta' \alpha
_0^{-1}$ and $\gamma = r \gamma' \alpha_0^{-1}$, we have
$\beta, \gamma \in (\mathbb {C}*H)_\sigma [X]$ and
\[
\beta \sum_i \alpha _i X^i + \gamma \sum_i \alpha_i
\zeta^i X^i = r.
\]
Set $\alpha' = \sum_i \alpha_i \zeta^i x^i$.
Applying the homomorphism $\theta$, we now have $\beta \alpha + \gamma
\alpha' = r$.
By Lemma~\ref{Lbilinear} there exists $y\in \mathcal {G}_n$ such that
$yhy^{-1} = h$ for all $h\in H$ and $yxy^{-1} = \zeta x$, and then we
have $y\alpha y^{-1} = \alpha'$.  Thus $\ker \alpha' = y (\ker \alpha
) y^{-1}$, consequently $\mathcal N(\alpha') = y \mathcal N(\alpha )
y^{-1}$ and using Lemma~\ref{Ltrace}\eqref{Ltrace21},
we deduce that $\tau \mathcal N(\alpha') = \tau N\mathcal
(\alpha) > \nu/2$.

Suppose $f \in \ker\alpha \cap \ker\alpha'$.  Then $\alpha f =
\alpha'f = 0$, hence $r f = 0$ because $r = \beta \alpha +
\gamma \alpha'$, and we can invoke
our inductive hypothesis to deduce that $f=0$.   Therefore $\ker\alpha
\cap \ker\alpha' = 0$.  If $\pi$ is the projection onto
$\overline{ \ker\alpha + \ker\alpha' }$, we now see from
Lemma~\ref{Ltrace}\eqref{Ltrace5} that
\[
\tau \pi = \tau \mathcal N(\alpha) + \tau \mathcal N(\alpha') >
\nu/2 + \nu/2 = \nu.
\]
Using Lemma~\ref{Lore2} we may choose $\delta$ so that
$0 \ne \delta \in \mathbb {C}*G \alpha
\cap \mathbb {C}*G \alpha'$, and
then $\ker\delta \supseteq \ker \alpha + \ker \alpha'$, hence $\tau
\mathcal N(\delta) \ge \tau \pi > \nu $ by
Lemma~\ref{Ltrace}\eqref{Ltrace4}.  This contradicts the definition of
$\nu$ and \eqref{Lmain1} is proven.

\eqref{Lmain2}
This follows from \eqref{Lmain1} and the remarks immediately preceding
Theorem~\ref{Tmain}.

\eqref{Lmain3}  This follows from \eqref{Lmain2},
Lemma~\ref{Lore2} and the comments
immediately preceding Lemma~\ref{Lore1}.
\end{proof}

\bibliographystyle{amsplain}
\bibliography{linind}

\providecommand{\bysame}{\leavevmode\hbox to3em{\hrulefill}\thinspace}
\begin{thebibliography}{1}

\bibitem{arveson}
W.~Arveson, \emph{An invitation to {C}*-algebra}, Graduate Texts in
  Mathematics, vol.~39, Springer-Verlag, Berlin-New York, 1976.

\bibitem{berberian}
S.~K. Berberian, \emph{Baer $*$-rings}, Grundlehren, vol. 195, Springer-Verlag,
  Berlin-New York, 1972.

\bibitem{ber}
\bysame, \emph{The maximal ring of quotients of a finite von {N}eumann
  algebra}, Rocky Mountain J. Math. \textbf{12} (1982), 149--164.

\bibitem{heilram}
C.~Heil, J.~Ramanathan, and P.~Topiwala, \emph{Linear independence of
  time-frequency translates}, Proc. Amer. Math. Soc. \textbf{124} (1996),
  2787--2795.

\bibitem{kadison}
R.~V. Kadison and J.~R. Ringrose, \emph{Fundamentals of the theory of operator
  algebras, volume 1, elementary theory}, Pure and Applied Mathematics Series,
  vol. 100, Academic Press, London-New York, 1983.

\bibitem{zero}
P.~A. Linnell, \emph{Zero divisors and group von {N}eumann algebras}, Pacific
  J. Math. \textbf{149} (1991), 349--363.

\bibitem{rieffel}
Marc~A. Rieffel, \emph{Von {N}eumann algebras associated with pairs of lattices
  in {L}ie groups}, Math. Ann. \textbf{257} (1981), 403--418.

\bibitem{stein}
E.~M. Stein, \emph{Harmonic analysis: real-variable methods, orthogonality, and
  oscillatory integrals}, Princeton Mathematical Series, vol.~43, Princeton
  University Press, Princeton, N.J., 1993.

\end{thebibliography}

\end{document}